\documentclass[11pt]{amsart}

\usepackage[all]{xypic}

\usepackage{latexsym}

\usepackage{amssymb}
\usepackage{amsfonts}
\usepackage{amscd}
\usepackage{amsmath,amsthm}

\newtheorem{lemma}{Lemma}[section]

\newtheorem{lem}[lemma]{Lemma}
\newtheorem{prop}[lemma]{Proposition}
\newtheorem{thm}[lemma]{Theorem}
\newtheorem{cor}[lemma]{Corollary}

{

}

\theoremstyle{definition}

\theoremstyle{remark}

\numberwithin{equation}{section}

\newenvironment{pf}{\noindent{\bf Proof.}}{\hfill $\square$\medskip}

\def\CC{{\mathbb C}}

\def\PP{{\mathbb P}}

\def\RR{{\mathbb R}}

\def\ZZ{{\mathbb Z}}

\def\0ol{{\bar 0}}
\def\1ol{{\bar 1}}
\def\2ol{{\bar 2}}
\def\ol2{{\bar 2}}
\def\3ol{{\bar 3}}
\def\4ol{{\bar 4}}
\def\5ol{{\bar 5}}
\def\6ol{{\bar 6}}
\def\7ol{{\bar 7}}
\def\8ol{{\bar 8}}
\def\9ol{{\bar 9}}

\def\bold0{{\bf 0}}
\def\bold1{{\bf 1}}
\def\bold2{{\bf 2}} 
\def\bold3{{\bf  3}}
\def\bold4{{\bf 4}}
\def\bold5{{\bf 5}}
\def\bold6{{\bf 6}}
\def\bold7{{\bf 7}}
\def\bold8{{\bf 8}}
\def\bold9{{\bf 9}}

\def\P2Skly{\PP^2_{Skly}}

\def\End{\operatorname {End}}

\def\GL{\operatorname {GL}}

\def\PGL{\operatorname {PGL}}

\def\th{\operatorname {th}}    

\def\Aut{\operatorname{Aut}}

\def\dim{\operatorname{dim}}

\def\End{\operatorname{End}}

\def\gcd{\operatorname{gcd}}
\def\GKdim{\operatorname{GKdim}}

\def\id{\operatorname{id}}

\def\liminj{\varinjlim}

\def\NS{{\sf NS}}

\def\Pic{\operatorname{Pic}}

\def\Spec{\operatorname{Spec}}

\def\sup{\operatorname{sup}}

\def\trdeg{\operatorname{trdeg}}

\def\ul1{\operatorname{\underline{1}}}

\def\l{\leftarrow}

\def\a{\alpha}
\def\b{\beta}

\def\l{\lambda}

\def\s{\sigma}

\def\ve{\varepsilon}

\def\L{\Lambda}

\def\fm{{\mathfrak m}}

\def\cal{\mathcal}

\def\cI{{\cal I}}

\def\cL{{\cal L}}

\def\cO{{\cal O}}

\def\cW{{\cal W}}

\def\cZ{{\cal Z}}

\def\dirlim{\mathop{\vtop{\baselineskip -100pt\lineskip -1pt\lineskiplimit 0pt
\setbox0\hbox{lim}\copy0\hbox to \wd0{\rightarrowfill}}}\limits}
\def\invlim{\mathop{\vtop{\baselineskip -100pt\lineskip -1pt\lineskiplimit 0pt
\setbox0\hbox{lim}\copy0\hbox to \wd0{\leftarrowfill}}}\limits}

\def\I11{{1 \kern -0.8pt \! \mbox{l}}}
\def\mumu{{\mu\kern-4.2pt\mu}}
\def\bfmu{{\mu\kern-4.2pt\mu}}
\def\2slash{\backslash \! \backslash}

\def\boxtimes{\setbox0\hbox{$\Box$}\copy0\kern-\wd0\hbox{$\times$}}

\begin{document}

\title[Free algebras and automorphisms of surfaces]{Free algebras arising from positive-entropy
automorphisms of  surfaces
}

\author{S. Paul Smith}
\address{ Department of Mathematics, Box 354350, Univ.  Washington, Seattle, WA 98195}
\email{smith@math.washington.edu}



\maketitle

\begin{abstract}
Let $X$ be a smooth projective 
surface defined over a field $k$ and let $k(X)$ be its field of rational functions.
Let $\s \in \Aut_k(X)$. This paper proves there is an integer $n \ge 0$ and elements $a,b \in k(X)$ such that 
the subalgebra of the skew Laurent extension $k(X)[t^{\pm 1};\s]$ generated by $at^n$ and $bt^n$ is 
a free algebra if and only if the spectral radius for the action of $\s^*$ on the N\'eron-Severi group of $X$ is $>1$.
Thus, when $\s$ is an automorphism of a smooth complex projective surface $X$, 
$\CC(X)[t^{\pm 1};\s]$ has a free subalgebra on $\ge 2$ variables if and only if the topological entropy of $\s$ is positive.  
Furthermore,  if $\s$ is an automorphism of $k(X)$, then $k(X)[t^{\pm 1};\s]$ contains a 
free subalgebra if and only if the dynamical degree of $\s$ is $>1$; in this situation, $\s$ might not be induced by
an automorphism of any smooth projective surface. These results are used to show that certain twisted 
homogeneous coordinate rings of smooth projective surfaces contain free subalgebras. 
\end{abstract}

\section{Introduction}

Throughout, $k$ is an uncountable algebraically closed field; $X$ is a smooth projective  surface over $k$;
$\s:X \dashrightarrow X$ is a birational map; $k(X)$ is the field of rational functions on $X$; and 
$\s$ also denotes the automorphism of $k(X)$ defined by $\s(f)=f\circ \s$.   

Whenever we say ``free algebra'' we mean ``free algebra on $\ge 2$ generators''.

\subsection{}



The {\sf skew Laurent polynomial ring} $k(X)[t^{\pm 1};\s]$  is defined to be the vector space 
$k(X) \otimes_k k[t^{\pm 1}]$ with multiplication being the bilinear extension of $(f \otimes t^m)(g \otimes t^n)=f\s^m(g)\otimes t^{n+m}$ for $f,g \in k(X)$.  We make 
$k(X)[t^{\pm 1};\s]$  a $\ZZ$-graded algebra by placing $k(X)$ in degree zero and setting
$\deg(t)=1$. 

If $u$ and $v$ are elements of $k(X)[t^{\pm 1};\s]$ we write $k\{u,v\}$ for the subalgebra they generate.

\subsection{}
We write $\NS(X)$ for the N\'eron-Severi group of $X$ and define $\NS(X)_\RR:=\NS(X) \otimes_\ZZ \RR$. 
There is a finite set $S \subset X$ such that the restriction of $\s$ to $U:=X-S$ is a morphism $\s|_{{}_U}:U \to X$. If $D$ is a 
divisor on $X$, the Zariski closure of $(\s|_{{}_U})^*D$ is a divisor on $X$ that we denote by $\s^*(D)$.  Thus
$\s$ induces a well-defined linear automorphism $\s^*:\NS(X)_\RR \to \NS(X)_\RR$.

We will always write $\rho(\s^*)$ for the spectral radius of $\s^*$.

The intersection pairing on $\Pic(X)$ induces a non-degenerate symmetric bilinear form on $\NS(X)_\RR$ 
that we denote by $(u,v) \mapsto u \cdot v$. The Hodge Index Theorem tells us that the signature of 
the corresponding quadratic form is $(1,d-1)$ where $d=\dim \NS(X)_\RR$. 

\subsection{}

Suppose $\s$ is a morphism, not just a birational map.
Then  $\s^*$ is an orthogonal transformation with respect to the quadratic form on $\NS(X)_\RR$. 
Furthermore, $\rho(\s^*) \ge 1$ and, if $\rho(\s^*)>1$, then 
$\rho(\s^*)$ and $\rho(\s^*)^{-1}$ are simple eigenvalues of $\s^*$.
 
These properties need not hold if $\s$ is just a birational map. For example, if $\s: \PP^2 \dashrightarrow  \PP^2$ is 
the map $\s(x,y,z)=(x^{-1},y^{-1},z^{-1})$, then $\s^*$ acts as multiplication by 2 on $\NS(\PP^2) \cong \ZZ$. 
 
 \subsection{}
 \label{sect.main.result}
 The first new result in this paper is part (1) of the following theorem. 
 Part (2) of it follows easily from  
 results of Artin-Van den Bergh \cite{AV}, Keeler \cite{K}, and Rogalski-Zhang \cite{RZ}.

\begin{thm}
\label{thm.main}
Let $X$ be a smooth projective surface over an uncountable algebraically closed field $k$. Let $\s \in \Aut(X)$. 
\begin{enumerate}
  \item 
  Suppose $\rho(\s^*) >1$. If $n$ is such that $\rho(\s^*)^n \ge 5+2\sqrt{6}$, then $k\{at^n,bt^n\}$ is a free subalgebra of $k(X)[t^{\pm 1};\s]$ for all $a,b \in k(X)$ such that the divisor of zeroes of $ab^{-1}$ is very ample.  
  \item 
If $\rho(\s^*) =1$, then $k(X)[t^{\pm 1};\s]$ does not contain a free subalgebra on $\ge 2$ variables.
\end{enumerate} 
\end{thm} 

The requirement that $\rho(\s^*)^n \ge 5+2\sqrt{6}$ can be replaced by $\rho(\s^*)^n \ge 2+\sqrt{3}$ when $\dim\NS(X)_\RR=2$.
See Corollary \ref{cor.improve} for a stronger statement.

\subsubsection{Sometimes the hypothesis that $k$ is uncountable and algebraically closed can be removed.}
Suppose that $X$ is a smooth projective surface over an arbitrary field $F$. Let $k$ be an uncountable algebraically closed
field that contains a copy of $F$. Let $X_k$ be the surface obtained by base extension from $\Spec(F)$ to $\Spec(k)$. Suppose $X_k$ is smooth and irreducible.  If $a,b \in F(X)$ and 
the divisor of zeroes of $ab^{-1}$ viewed as an element of $k(X)$ is very ample, then $F\{at^n,bt^n\}$ is a free subalgebra of 
$F(X)[t^{\pm 1};\s]$ if $n$ is such that $\rho(\s^*)^n \ge 5+2\sqrt{6}$.
  
 Similar remarks apply to the other main result in this paper, namely Theorem \ref{thm.2}.  
 The proofs of Theorems \ref{thm.main} and \ref{thm.2} make essential use of the hypothesis that $k$ is uncountable and algebraically closed, however, it is obvious that if $a$ and $b$ are in $F(X)$ and the subalgebra $k\{at^n,bt^n\}$ of 
 $k(X)[t^{\pm 1};\s]$ is free so is the subalgebra $F\{at^n,bt^n\}$ of 
 $F(X)[t^{\pm 1};\s]$.

\subsection{}

 Theorem \ref{thm.main} provides information about the twisted homogenous coordinate ring  $B(X,\cL,\s)$ associated to a very ample invertible $\cO_X$-module $\cL$.  Such rings play a central role in non-commutative algebraic geometry. 
They are non-commutative analogues of the commutative rings $\oplus_{n \ge 0}  H^0(X,\cL^{\otimes n})$.
If $\cL \cong \cO_X(D)$, then $B(X,\cL,\s)$ is isomorphic to the subalgebra  
$$
\bigoplus_{n=0}^\infty H^0\big(X,\cO_X(D+\s^{-1}D+\cdots+\s^{1-n}D)\big)t^n \; \subset \; k(X)[t^{\pm 1};\s].
$$

Part (1) of the following corollary follows from  Theorem \ref{thm.main}(1);
part (2) is due to Artin and Van den Bergh \cite{AV}. 

\begin{cor}
\label{cor.main}
Let $X$ be a smooth projective surface over an uncountable algebraically closed field $k$, and $\s \in \Aut(X)$. 
Let $\cL$ be a very ample invertible $\cO_X$-module. 
\begin{enumerate}
  \item 
If $\rho(\s^*) >1$,  then $B(X,\s,\cL)$ has a free subalgebra on two homogeneous elements of degree $n$ where $n$ has
the properties in Theorem \ref{thm.main}. 
  \item 
If $\rho(\s^*) =1$, then $B(X,\s,\cL)$ is a finitely generated, left and right noetherian, domain having finite 
Gelfand-Kirillov dimension.
\end{enumerate} 
\end{cor}

\subsection{}
An essential step in the proof of Theorem \ref{thm.main}(1) is proving the following: if $H\subset X$ is an 
effective   very ample divisor and $n$ is as in Theorem \ref{thm.main}(1), then there is a smooth irreducible curve 
$C$ such that  the intersection numbers satisfy the inequality 
$\s^{n(j+1)}(H)\cdot C\, > \, 2 \s^{nj}(H\cdot C)$ for all $j \ge 0$.  

There are many results in the literature about the behavior of the intersection numbers
$\s^i(H) \cdot C$ as $i \to \infty$, but those estimates do not seem sufficient to prove Theorem \ref{thm.main};
our proof requires the strict inequality in the previous paragraph.

The fact that the iterated images of $H$ under the powers of $\s$ intersect $C$ at an exponential rate is a 
statement about the dynamical behavior of $\s$ on the curves lying on $X$. When $k=\CC$, 
Theorem \ref{thm.main} can be stated in terms of the entropy of $\s$ (see Theorem \ref{thm.ent}).

\subsection{Automorphisms of complex surfaces}

 \begin{thm}
[Gromov \cite{Grom03}, Yomdin \cite{Y}]
If $\s$ is an automorphism of a smooth complex projective surface $X$, then the topological entropy of 
$\s$ is $\log(\l)$  where $\l=\rho(\s^*)$.
\end{thm}


 \begin{thm}
 \label{thm.ent}
Let $X$ be a smooth  complex projective surface and $\s \in \Aut(X)$. 
If the topological entropy of $\s$ is positive there is an integer $n$ and functions $a,b \in \CC(X)$ such that $\CC\{at^n,bt^n\}$ 
is a free subalgebra of $\CC(X)[t^{\pm 1};\s]$. If the topological entropy of $\s$ is zero, then $\CC(X)[t^{\pm 1};\s]$
does not contain a free subalgebra on $\ge 2$ variables.
\end{thm}

\subsection{}
As mentioned in \S\ref{sect.main.result}, Theorem \ref{thm.main}(2) follows rather easily from earlier results.
Those results, which are described in \S\ref{sect.proof}, play a central  role in non-commutative algebraic geometry.
Among other things, they imply that  every finitely generated 
subalgebra of $k(X)[t^{\pm 1};\s]$ has finite Gelfand-Kirillov dimension 
(GK-dimension) when $\rho(\s^*) =1$; however, a free algebra on $\ge 2$ variables has infinite GK-dimension so 
cannot be a subalgebra of  $k(X)[t^{\pm 1};\s]$ when $\rho(\s^*)=1$.

\subsection{Birational self-maps}
A natural algebraic question presents itself:  if $K/k$ is a  finitely generated field extension such that $\trdeg_k(K)=2$,
 and $\s \in \Aut_k(K)$, when does  $K[t^{\pm 1}; \s]$ contain a free subalgebra?
 
When $\s$ is induced by  an automorphism of a smooth projective surface $X$ such that $k(X)=K$, the question is answered by Theorem \ref{thm.main}. However, there are automorphisms of $k(\PP^2)$ that are not induced by an
automorphism of any smooth projective surface $X$ with the property that $k(X) \cong k(\PP^2)$. 

This leads us to the second part of this paper where a more precise version of part (1) of the following result is proved.

\begin{thm}
\label{thm.2}
Let $X$ be a smooth projective surface over an uncountable algebraically closed field $k$. 
Let $\s: X \dashrightarrow X$ be a birational map and define $\l(\s):=\lim_{n \to \infty} \rho\big((\sigma^n)^*\big)^{1/n}$.   
\begin{enumerate}
\item{}
 Suppose $\l(\s) >1$. If $n$ is such that $\l(\s)^n \ge 5+2\sqrt{6}$, then there are non-zero elements $a,b \in k(X)$
 such that $k\{at^n,bt^n\}$ is a free subalgebra of $k(X)[t^{\pm 1};\s]$.  
\item 
If $\l(\s)=1$, then $k(X)[t^{\pm 1};\s]$ does not contain a free subalgebra on $\ge 2$ variables.
\end{enumerate}
\end{thm}

Part (2) of the above theorem is an immediate consequence of the work of others.

\subsection{}
Section \ref{sect.egs} gives explicit examples of free subalgebras of $k(X)[t^{\pm 1};\s]$  for various K3 surfaces $X$.

Section \ref{ssect.PP2} exhibits some automorphisms $\s$ of $k(\PP^2)$ that, although they are not induced by 
an automorphism of any model of $k(\PP^2)$, have the property that $k(\PP^2)[t^{\pm 1};\s]$ contains a free subalgebra.  
We refer the reader to \cite{SPS} for more precise information about free subalgebras of $k(\PP^2)[t^{\pm 1};\s]$.
The methods in \cite{SPS}, which mostly involve valuations on $k(\PP^2)$,  are very different from those in this paper.

\subsection{Acknowledgements}
My interest in free subalgebras of $k(X)[t^{\pm 1};\s]$ was roused by 
Jason Bell, Kenneth Chan, and Dan Rogalski. I thank them all for useful and stimulating conversations.

This work was supported by the National Science Foundation under Grant No. 0932078 000 while the author was in residence at the Mathematical Science Research Institute (MSRI) in Berkeley, California, during the Noncommutative Algebraic Geometry and Representation Theory Program in the winter semester of 2013.

\section{The proof of Theorem \ref{thm.main}}
\label{sect.proof}

Subsections \ref{ssect.twhcr} and \ref{ssect.RZ} recall the results and definitions needed to prove 
Theorem \ref{thm.main}(2):  $k(X)[t^{\pm 1};\s]$  does not contain a free subalgebra on $\ge 2$ variables when 
$\rho(\s^*)=1$.  

In \S\ref{sect.main.thm} we prove Theorem \ref{thm.main}(1).

\subsection{Twisted homogeneous coordinate rings}
\label{ssect.twhcr}
Let $X$ be a projective scheme over an arbitrary field $k$, $\s \in \Aut(X)$, and $\cL$ an ample invertible $\cO_X$-module.
 The ring 
 $$
 B(X,\cL,\s) := \bigoplus_{n=0}^\infty H^0\big(X,\cL \otimes \s^*\cL \otimes \cdots \otimes (\s^*)^{n-1}\cL\big),
 $$
with multiplication defined by $a\cdot b:=a\otimes (\s^*)^m(b)$ when $a \in B_m$ and $b \in B_n$,
is called a {\sf twisted homogeneous coordinate ring} for $X$. The ring $B(X,\cL,\s)$ was introduced by Artin, Tate, and 
Van den Bergh \cite{ATV}. 

When $X$ is integral, $B(X,\cL,\s)$ is isomorphic to a graded subalgebra of $k(X)[t^{\pm 1};\s]$. 
If $D$ is a divisor such that $\cL \cong \cO_X(D)$, then
$$
B(X,\cL,\s) \cong \bigoplus_{n=0}^\infty  H^0\big(X,\cO_X(D+\s^{-1} D + \cdots + \s^{1-n}D)\big) t^n \, \subset \, k(X)[t^{\pm 1};\s].
$$
 
\subsection{}
In \cite{AV},  Artin and Van den Bergh  laid out the basic properties of twisted homogeneous coordinate rings  and investigated the case of surfaces in detail. The main result in their paper is that if $\cL$ is $\s$-ample, then
 a suitable quotient of the category of graded $B(X,\cL,\s)$-modules is equivalent to the category of quasi-coherent
  $\cO_X$-modules.  We will not make use of that result.

Although we will use the implications of Theorems \ref{thm.AV} and \ref{thm.K} the reader of this paper 
does not need to know the meaning of ``$\s$-ample''. Information about 
the definition, existence, and consequences of $\s$-ampleness can be found  in \cite{AV} and \cite{K}.

\begin{thm}
[Artin-Van den Bergh]
\cite[Thm. 1.7]{AV}
\label{thm.AV}
Let $\s$ be an automorphism of a smooth proper algebraic surface $X$.
If there is a $\s$-ample invertible $\cO_X$-module, then $\rho(\s^*)=1$.
\end{thm}

The {\sf Gelfand-Kirillov dimension}, or GK-dimension, of a $k$-algebra $A$ is
$$
\GKdim( A) := \sup_V \overline{\lim_{n \to \infty}} \log_n\!\big(\!\dim_k V^n\big)
$$
where the supremum is taken over all finite-dimensional subspaces $V\subset A$.
 
\begin{thm}
[Keeler]
\cite[Thm. 1.4]{K}
\label{thm.K}
Let $X$ be a projective scheme over an algebraically closed field, and let $\s$ be an automorphism of $X$. 
The following are equivalent: 
\begin{enumerate}
  \item 
 $\rho(\s^*)=1$;
  \item 
 there is a $\s$-ample invertible $\cO_X$-module;
 \item{}
 every  ample invertible $\cO_X$-module is $\s$-ample;
  \item 
  $B(X,\cL,\s)$ has finite GK-dimension for all ample invertible $\cO_X$-modules $\cL$;
   \item 
  $B(X,\cL,\s)$ is left and right noetherian for all ample invertible $\cO_X$-modules $\cL$.
\end{enumerate}  
\end{thm}
 
\subsection{}
\label{ssect.RZ}
In \cite{RZ}, Rogalski and Zhang  showed twisted homogeneous coordinate rings are ubiquitous by
proving that for a large class of graded rings, $A$, there is a canonical graded ring homomorphism
$A \to B(X,\cL,\s)$ that is surjective in large degree where $X$ is, as in \cite{ATV}, the solution to a moduli problem
involving a class of graded $A$-modules,
$\s$ is induced by the Serre twist (degree shift) on graded $A$-modules, and $\cL$ is $\s$-ample.

\begin{thm}
[Rogalski-Zhang]
\cite[Thm. 1.6]{RZ}
\label{thm.RZ}
Let $X$ be an integral projective scheme and $\s$ an automorphism of $X$. If $\rho(\s^*)=1$, then 
$$
\GKdim\!\big(k(X)[t^{\pm 1};\s]\big)\le \GKdim\!\big(B(X,\cL,\s)\big)+\dim(X) < \infty.
$$
In particular, if $\rho(\s^*)=1$, then $k(X)[t^{\pm 1};\s]$ does not contain a free algebra on $\ge 2$ variables.
\end{thm}

\subsection{}

Proofs of the following result can be found in \cite[\S5]{AV}, \cite[Thm. 5.1]{DF}, and \cite[\S3]{M1}.
 
\begin{prop}
\label{prop.q.form}
Let $(V,q)$ be a $d$-dimensional vector space over $\RR$ endowed with a quadratic form $q$ having signature $(1,d-1)$.
Let $\theta$ be an orthogonal automorphism of $V$ that preserves a $\ZZ$-lattice.  Let $\l$ denote the spectral radius of 
$\theta$. If $\l>1$, then 
\begin{enumerate}
\item
$\l$ and $\l^{-1}$ are simple eigenvalues for $\theta$; 
  \item
if $v$ is a $\l$- or $\l^{-1}$-eigenvector for $\theta$, then $q(v)=0$;
\item{}
every totally isotropic subspace of $(V,q)$ has dimension $1$. 
\end{enumerate}
\end{prop}

Let $\s$ be an automorphism of a smooth projective surface $X$. 
   By the Hodge Index Theorem, the signature of the quadratic form on $\NS(X)_\RR$ induced by the intersection pairing 
 is $(1,d-1)$ where $d=\dim \NS(X)_\RR$. 
 The action of $\s^*$ on $\NS(X)_\RR$ is orthogonal with respect to this form.
 We will apply Proposition \ref{prop.q.form} to the action of $\s^*$ on $\NS(X)_\RR$.


\subsection{}
\label{sect.main.thm}

In this section we prove Theorem \ref{thm.main}(1). 

First, we need Lemma \ref{lem.known}, which is probably well-known, Lemma \ref{lem.R.6.3} which 
is part of  the proof of \cite[Lemma 6.3]{Rog}, and Lemma \ref{lem.R.6.3-2nd} which is a simple consequence of 
Lemma \ref{lem.R.6.3}.

\begin{lem}
\label{lem.known}
Let $\s$ be an automorphism of a commutative ring $R$. Let $a,b,c,d \in R$ and consider the subalgebras $k\{at,bt\}$ and $k\{ct,dt\}$ of $R[t^{\pm 1};\s]$. If  $a,b,c,d$ are regular elements and $ad=bc$, then there is an algebra isomorphism
 $\Psi:k\{at,bt\} \cong k\{ct,dt\}$ such that $\Psi(at)=ct$ and $\Psi(bt)=dt$. In particular, if one of the algebras $k\{at,bt\}$, $k\{t,ab^{-1}t\}$, $k\{t,a^{-1}bt\}$, is free, then all three are.
\end{lem}
\begin{pf}
 Define $v=a^{-1}c=b^{-1}d$ and $v_m:=v\s(v)\cdots \s^{m-1}(v)$ for all $m \ge 1$. 

Let $\Psi:k\{at,bt\} \to k\{ct,dt\}$ be the $k$-linear extension of the map $\Psi(gt^m):=v_mgt^m$ for $g \in R$.
This is a $k$-algebra homomorphism because
\begin{align*}
\Psi(rt^m)\Psi(st^n) &  \, = \, v_m \s^m(v_n)r\s^m(s)t^{m+n}
\\
& \, = \, v_{m+n}r\s^m(s)t^{m+n} 
\\
& \, =\, \Psi(rt^mst^n).
\end{align*}
Note that $\Psi(at)=ct$ and $\Psi(bt)=dt$. 
There is a similar homomorphism $\Phi: k\{ct,dt\} \to k\{at,bt\}$ defined by $\Phi(gt^m):=v_m^{-1}gt^m$ for $g \in R$.
It is easy to check that $\Phi$ and $\Psi$ are mutual inverses.
\end{pf}

\subsubsection{Notation}
If $G$ is an abelian group we write $G^\times$ for the complement of the identity element. 

If $f$ is a non-constant rational function on a a smooth projective curve $C$, we write $\deg(f)$ for the degree of the 
rational map $f:C \dashrightarrow \PP^1$.

\begin{lem}
[Rogalski]
\cite[Lemma 6.3]{Rog}
\label{lem.R.6.3}
Let $C$ be a smooth projective curve over an uncountable 
algebraically closed field $k$. Let $T$ be a finite dimensional subspace of $k(C)$.
Let $U$ be a 2-dimensional subspace of $k(C)$ that contains $k$. If there is an element $f \in U^\times$ such that $\deg(f) >
\deg(g)$ for all $g \in T^\times$, then $\dim(TU)=2\dim(T)$.
\end{lem}
\begin{pf}
Suppose $f$ is such a function.
Since the set
$$
S:=\{p \in C \; | \; \hbox{some $g \in T$ has a pole at $p$}\}
$$
 is finite, there is $\a \in k$ such that $f+\a$ does not
vanish at any point of $S$.   Since $\deg(f)=\deg(f+\a)$ we can, and will, replace $f$ by $f+\a$. Let $g \in T^\times$. 
Since $f$ does not vanish at any point where $g$ has a pole, $fg$ has a zero everywhere $f$ does.
Therefore $\deg(fg)\ge \deg(f)>\deg(g)$.
Hence $fg \notin T$. Thus $T \cap fT=0$ and $UT=T+fT = T \oplus fT$.
\end{pf}

\begin{lem}
\label{lem.R.6.3-2nd}
Let $C$ be a smooth projective curve over an uncountable algebraically closed field $k$. Let $U_i=k\oplus kf_i$, $0 \le i \le n$, be 2-dimensional subspaces of $k(C)$.  
If  $\deg(f_{i+1})\ge 2\deg(f_i)$ for all $i \ge 0$, then
$$
\dim_k(U_0U_1 \cdots U_n)=2^{n+1}
$$
\end{lem}
\begin{pf}
The lemma is certainly true for $n=0$. We now argue by induction using 
Lemma \ref{lem.R.6.3} with $U=U_n$ and $T=U_0\cdots U_{n-1}$.
Suppose  $\dim_k(T)=2^n$.
The hypothesis on the degrees of the $f_i$s implies that $\deg(f_i)>2 \deg(g)$ for all non-zero $g \in U_{i-1}$.
Therefore
$$
\deg(f_n)\,>\,\deg(f_{n-1}) + \cdots + \deg(f_0) \ge \deg(g)
$$
for all non-zero $g \in U_0\cdots U_{n-1}$. The result now follows from Lemma \ref{lem.R.6.3}. 
\end{pf}

\begin{thm}
\label{main.thm}
Let $X$ be a smooth projective surface over an uncountable algebraically closed field $k$. 
Let $H$  be an effective very ample  divisor on $X$ and choose $h \in k(X)$  such that $(h)_0=H$.  
If $\s$ is an automorphism of $X$ such that $\rho(\s^*)^n \ge 5+2\sqrt{6}$,
then 
\begin{enumerate}
  \item 
$k\{t^n,ht^n\}$ is a free subalgebra of $k(X)[t^{\pm 1};\s]$;
  \item 
  $k\{at^n,bt^n\}$ is a free subalgebra of $k(X)[t^{\pm 1};\s]$ if $a,b \in k(X)$ are such that $a^{-1} b=h$. 
\end{enumerate}
\end{thm}
\begin{pf} 
If (1) is true, then Lemma \ref{lem.known} implies (2) is. We will prove (1).

Because the subalgebra of $k(X)[t^{\pm 1};\s]$ generated by $t^{\pm n}$ is isomorphic to $k(X)[t^{\pm 1};\s^n]$ we can, and 
will, replace $\s$ by $\s^n$ where $n$ is chosen as in the statement of the theorem. Thus, we will prove that  $k\{t,ht\}$ is a free algebra under the assumption that $\rho(\s^*)=\l \ge 5+2 \sqrt{6}$.

The algebra $k\{t,ht\}$ is graded by placing $k(X)$ in degree zero and $t$ in degree one. The degree-one component
of $k\{t,ht\}$ is $Vt$ where $V=k+kh$. The degree-$(j+1)$ component of $k\{t,ht\}$ is 
$$
(Vt)^{j+1}= V \s(V) \cdots \s^n(V)t^j.
$$
We will prove that $k\{t,ht\}$ is free by showing that
$$
\dim_k \! \big(V \s(V) \cdots \s^j(V)\big) = 2^{j+1}
$$
for all $n \ge 0$. 

Let $d=\dim\NS(X)_\RR$. Because $\rho(\s^*)>1$, $\s^*$ has at least two different eigenvalues. Hence $d \ge 2$. 

 Let $e_{{}_+},e_{{}_-} \in \NS(X)_\RR$ be eigenvectors for $\s^*$ with $\s^*(e_\pm)=\l^{\pm 1}e_{\pm}$. 
 By Proposition \ref{prop.q.form}, $ e_{{}_+}\!\! \cdot e_{{}_+} = e_{{}_-}\! \! \cdot e_{{}_-}=0$ and, because
$\RR e_{{}_+} \oplus \RR e_{{}_-}$ is not totally isotropic, $e_{{}_+}\!\! \cdot  e_{{}_-} \ne 0$.  Let $W=(\RR e_{{}_+} \oplus \RR e_{{}_-})^\perp$.
Then $\NS(X)_\RR =  \RR e_{{}_+} \oplus \RR e_{{}_-} \oplus W$. 
Since $\s^*$ is orthogonal,  $\s^*(W)=W$. By the Hodge Index Theorem, the intersection form is negative definite on $W$. 

Write $[H]= \a e_{{}_+} + \b e_{{}_-} + w$ for the class of $H$ in $\NS(X)_\RR$ where 
$w \in W$. Because $H\!\! \cdot \!  H >0$ and $w\cdot w\le 0$, $\a\b(e_{{}_+}\!\! \cdot \!  e_{{}_-})> 0$. 
Thus we can, and will, replace $e_{{}_+}$ and $e_{{}_-}$ by multiples of themselves  so that $[H]= e_{{}_+} \, + \, e_{{}_-} \, +\, w$. With this notation,
\begin{align}
\label{intersect.no}
 (\s^{*})^j[H]\!  \cdot \!  [ H]
&= (\l^j e_{{}_+} +  \l^{-j} e_{{}_-} + (\s^*)^jw ) \! \cdot \!  (e_{{}_+} + e_{{}_-} + w)
\notag
\\
&=(\l^j + \l^{-j})(e_{{}_+}\!\! \cdot \!  e_{{}_-}) +  (\s^*)^jw \! \cdot \!  w
\end{align}

Because $H \cdot \!  H>0$,  $2(e_{{}_+}\!\! \cdot \!  e_{{}_-})+ w\! \cdot \!  w > 0$. Therefore
$$
-2 < \bigg(\frac{w \cdot w}{e_{{}_+}\!\! \cdot \!  e_{{}_-}  }\bigg) \le 0.
$$
The Cauchy-Schwarz inequality in $\NS(X)_\RR$ gives
$$
\Big((\s^*)^j w \cdot w \Big)^2 \, \le \, \Big((\s^*)^j w \cdot (\s^*)^j w \Big) (w \cdot w) \, = \, (w \cdot w) ^2
$$
for all $j \ge 0$. 
Therefore
$$
\bigg(\frac{1}{e_{{}_+}\!\! \cdot \!  e_{{}_-}  }\bigg) \Big(2(\s^*)^j w \cdot w - (\s^*)^{j+1} w \cdot w \Big) \, \le \,
 -3 \bigg(\frac{w \cdot w}{e_{{}_+}\!\! \cdot \!  e_{{}_-}  }\bigg)  < 6.
$$
The hypothesis on $\l$ implies $(\l-2)+\l^{-1}(1-2\l) > 6 $. Since $1-2\l$ is negative it follows that 
$$
\l^j(\l-2)+\l^{-j-1}(1-2\l) \,  > \, \bigg(\frac{1}{e_{{}_+}\!\! \cdot \!  e_{{}_-}  }\bigg) \Big(2(\s^*)^j w \cdot w - (\s^*)^{j+1} w \cdot w \Big).
$$
  Rearranging the previous equation gives 
\begin{equation}
\label{inequal.5}
  \l^{j+1}+\l^{-j-1}  +  \frac{(\s^*)^{j+1}w\! \cdot \!  w}{e_{{}_+}\!\! \cdot \!  e_{{}_-}}  \, > \, 2\Bigg(\l^j+\l^{-j} + \frac{(\s^*)^{j}w\! \cdot \!  w}{e_{{}_+}\!\! \cdot \!  e_{{}_-}} \Bigg).
 \end{equation}
Combining (\ref{inequal.5})  and (\ref{intersect.no}) gives
$$
 (\s^{*})^{j+1}[H]\!  \cdot \!  [ H] \, > \,  2(\s^{*})^j[H]\!  \cdot \!  [ H]
 $$
 for all $j \ge 0$.
This inequality remains true if $H$ is replaced by a divisor that is linearly equivalent to $H$. 

Because $H$ is very ample and effective, there is $h \in k(X)$ whose divisor of zeroes is $H$. We now consider $V =k+kh$
and $\s^j(V)=k+k\s^j(h)$. 

It follows from  Bertini's Theorem \cite[8.20.2]{Ha} and the uncountability of $k$ that there is a smooth irreducible curve 
$C \in |H|$  such that  none of the  functions  $\sigma^j(h^{\pm 1})$ vanishes on $C$. Let $C$ be such a curve.

Let $\fm_C$ denote the maximal ideal in the discrete valuation ring $\cO_{X,C}$. Thus, $\cO_{X,C}/\fm_C \cong k(C)$.
We will show that $\s^j(V) \subset \cO_{X,C}$  and $\fm_C \cap \s^j(V)=\{0\}$.

 
 Since $H$ is very ample,  $\s^j(H) \cdot C >0$ for all $j \in \ZZ$. Since $C$ and 
 $\s^j(H)$ do not have a common irreducible component, $\s^j(H) \cap C$ is a finite non-empty set. Hence $\s^j(h)$ vanishes
 at a finite non-empty subset of $C$. It follows that $\a+\s^j(h) \in \cO_{X,C} - \fm_C$ for all $\a \in k$ and all $j \in \ZZ$. 
 In other words, $\s^j(V) \subset \cO_{X,C}$ and its image in $k(C)$, which we will denote by 
$U_j$, has dimension 2.  

From now on we only consider the case $j \ge 0$. 
Because 
$$
\deg\big({\s^j}(h)\big\vert_{{}_C}\big) = \s^j(h) \cdot \ C \;
> \; 2(\s^{j-1}(h) \cdot  C)=2 \deg\big({\s^{j-1}}(h)\big\vert_{{}_C}\big),
$$
$\deg(f)>2\deg(g)$ for all non-zero $f \in U_{j}$, all non-zero $g \in U_{j-1}$, and all $j$. 
It now follows from Lemma \ref{lem.R.6.3-2nd} that $\dim_k(U_0U_1 \cdots U_n)=2^{n+1}$. 

But $U_j$ is the image of $ \s^j(V)$ in $k(C)$ so $U_0U_1 \cdots U_n$ is the image of $V\s(V)\cdots \s^n(V)$ in $k(C)$. 
We conclude that
$$
\dim_k\big(V\s(V)\cdots \s^n(V)\big)=2^{n+1}
$$
and therefore $k\{t,ht\}$ is a free algebra.
\end{pf}

Sometimes Theorem \ref{main.thm} does not give the smallest $n$. For example, the next result gives a smaller $n$
 if $\dim \NS(X)_\RR=2$ or the eigenvalues of $\s^*$ are $\rho(\s^*)^{\pm 1}$ and roots of unity. 

\begin{cor}
\label{cor.improve}
Adopt the notation in Theorem \ref{main.thm}. 
If there is an $n$  such that the eigenvalues of 
$(\s^*)^n$ belong to $\{1,\rho(\s^*)^n, \rho(\s^*)^{-n}\}$ and $\rho(\s^*)^n > 2+\sqrt{3}$, 
then $k\{at^n,bt^n\}$ is a free algebra if $(a^{-1}b)_0$ is very ample.
\end{cor}
\begin{pf}
Replace $\s$ by $\s^n$. The hypotheses imply that  if $W \ne 0$, the only eigenvalue for the action of $\s^*$ on $W$ is 1. 
As remarked in the proof of \cite[Lemma 5.4]{AV}, this implies that $\s^*$ acts as the identity on $W$. Therefore
(\ref{intersect.no}) can be replaced by the statement that 
$$
 (\s^{*})^n[H]\!  \cdot \!  [ H] =(\l^n + \l^{-n})(e_{{}_+}\!\! \cdot \!  e_{{}_-}) +  (w\! \cdot \!  w)
$$
for all $n \ge 0$.

The hypothesis that $\l > 2 + \sqrt{3}$ implies $\l(\l-2)+(1-2\l) > 0$. Thus, 
$$
\l^{2n+1}(\l -2) + (1-2\l) \, > \, 0 \, \ge \, \l^{n+1} \bigg(\frac{w\! \cdot \!  w}{e_{{}_+}\!\! \cdot \!  e_{{}_-}}\bigg)
$$ 
 for all $n \ge 0$. Rearranging this gives 
$$
  \l^{n+1}+\l^{-n-1}  +  \frac{w\! \cdot \!  w}{e_{{}_+}\!\! \cdot \!  e_{{}_-}}  \, > \, 2\Bigg(\l^n+\l^{-n} + \frac{w\! \cdot \!  w}{e_{{}_+}\!\! \cdot \!  e_{{}_-}} \Bigg).
  $$
Therefore
 $$
 ( \s^*)^{n+1}[H] \! \cdot \!  [H ] \, > \,  2 (\s^*)^{n}[H] \! \cdot \!  [H] 
 $$
  for all $n \ge 0$. 
  
 The rest of the proof proceeds as in Theorem \ref{main.thm}.
 \end{pf}

By \cite[Remark 2.6]{K1}, if $\l$ is a quadratic integer, then there are eigenvectors $e_{\pm}$ as in the proof of Theorem \ref{main.thm} with  the additional property that $e{{}_+}+e_{{}_-}$ is nef and big. 

\begin{cor}
\label{cor.vample}
Let $X$ be a smooth projective surface over an uncountable algebraically closed field $k$. Let $\s \in \Aut(X)$
and suppose  $\rho(\s^*)>1$. 
If $\cL$ is a very ample invertible $\cO_X$-module, then $B(X,\cL,\s)$ contains a graded free subalgebra.
\end{cor}
\begin{pf}
Let $D$ be an effective divisor such that $\cL \cong \cO_X(D)$. Then $B(X,\cL,\s)$ is isomorphic to  $B(X,\cO_X(D),\s)$
which, by its definition, is isomorphic to the graded subalgebra 
$$
B:=\bigoplus_{j=0}^\infty H^0\big(X,\cO_X(D+\s^{-1}D+\cdots+\s^{-j+1}D)\big)t^j
$$
of $k(X)[t^{\pm 1};\s]$. 

We will identify $B(X,\cO_X(D),\s)$ with this subalgebra of $k(X)[t^{\pm 1};\s]$. 

Choose $n$ such that the hypotheses of Theorem \ref{main.thm} are satisfied. Define $H:=D+\s^{-1}D+\cdots+\s^{-n+1}D$.
Since $H$ is very ample and effective, there is an element $h \in k(X)$ such that $(h)_\infty = H$. 
Since $(k+kh)t^n\subset B_n$, $B(X,\cO_X(D),\s)$ contains the subalgebra $k\{t^n,ht^n\}$ of $k(X)[t^{\pm 1};\s]$. 
By Theorem \ref{main.thm}, $k\{t^n,h^{-1}t^n\}$ is a free algebra. Since $k\{t^n,ht^n\} \cong k\{t^n,h^{-1}t^n\}$, $B(X,\cO_X(D),\s)$ contains a free subalgebra generated by two homogeneous elements of degree $n$. 
\end{pf} 

We have been unable to prove that the conclusion of 
Corollary \ref{cor.vample} holds when $\cL$ is ample, not just very ample.

\subsection{Additional remarks}


\begin{lem}
\cite[Lemma 3.2]{K}
\label{lem.DK}
If $\lambda$ is the spectral radius of $\s^*$, then  $\s^*$ has an eigenvector $v \in \NS(X)_\RR$ with eigenvalue 
$\lambda$ that belongs to the cone generated by
the numerically effective divisors.
\end{lem}

\begin{lem}
If $\s^*$ has an eigenvalue $>1$ it does not fix the class of any ample divisor in $NS(X)_\RR$.
\end{lem}
\begin{pf}
Suppose $\l>1$ is an eigenvalue of $\s^*$ and let $v$ be as in Lemma \ref{lem.DK}. 
If $D$ is an ample divisor such that $\s^*D=D$ in $NS(X)_\RR$, then
$\l(D \cdot v)=(D \cdot \s^*v)=((\s^*)^{-1}D \cdot v) =(D \cdot v)$. But  $v$ is in the interior of the effective cone so, by  
the Nakai-Moishezon criterion, $(D \cdot v)>0$.  This is absurd. 
\end{pf}

Let $K_X$ denote the canonical divisor. Every automorphism of $X$ fixes the class of $K_X$ in $NS(X)_\RR$ so,
by the previous lemma neither $K_X$ nor $-K_X$ can be ample when $\s^*$ has an eigenvalue $>1$. 

The next result is probably well-known. The proof is an easy generalization of \cite[Cor. 2.3]{S}.

\begin{prop}
Let $\s$ be an automorphism of a smooth projective surface $X$. Suppose there is $v\in \NS(X)_\RR$ 
that is contained in the ample cone and $\s^*v=\l v$ where $|\l| \ne 1$.  
Then every infinite $\s$-orbit in $X$ is Zariski-dense.
\end{prop}
\begin{pf}
Suppose the result is false. Then there is a point $p \in X$ such that $P:=\{ \s^m(p) \; | \; m \in \ZZ \}$ is infinite and its closure 
is the union of a finite number of irreducible curves, say $C_1 \cup \ldots \cup C_n$. 
For each $i$, the sets $P \cap C_i$, $\s(P) \cap C_i$, and $P \cap \s^{-1}(C_i)$ have the same cardinality. 
But $P$ meets each $C_i$ at infinitely many points so the closure of $\s^{-1}(C_i)$ is equal to some $C_j$.  
Hence $\s$ permutes the $C_i$s and therefore $\s^r(C_j)=C_j$ for some $r>0$ and some $C_j$. 

Since $v$ is ample,  $(v\cdot C_j)>0.$ But $\l^r(v \cdot C_j) = (\s^r(v) \cdot C_j)=
  (v \cdot \s^{-r}C_j) = (v \cdot C_j)$ so $(v \cdot C_j)=0$. These two contradictory facts imply the truth of the lemma. 
\end{pf}

Xie Junyi \cite[Thm. 1.1]{X} proves a much stronger result: if $\s:X \dashrightarrow X$ is a birational self-map of a smooth projective surface over an algebraically closed field whose characteristic is not 2 or 3, then the set of non-critical periodic  points
is Zariski dense in $X$. If, in addition, the base field is algebraically closed and $\l(\s)>1$,\footnote{The number $\l(\s)$ is defined in \S\ref{sect.dyn.deg}).} then there is a point $x \in X$ such that
$\s^n(x)$ is defined for all $n \in \ZZ$ and $\{\s^n(x) \; | \; n \in \ZZ\}$ is Zariski dense in $X$ \cite[Thm. 1.4]{X}.

\section{non-geometric birational maps}
\label{sect.birat.maps}

Throughout \S\ref{sect.birat.maps}, $k$ is an uncountable algebraically closed
field, $K/k$ is a finitely generated field extension of transcendence degree 2, $\s \in \Aut_k(K)$, 
and $X$ and $Y$  are smooth projective surfaces over $k$.

\subsection{}
This section addresses the question of when $K[t^{\pm 1};\s]$ contains a free subalgebra. 
If $\s$ is induced by an automorphism of a surface whose function field is $K$, then Theorem \ref{main.thm} 
answers the question. However, not all automorphisms of all $K$ are so induced. 

The main result in 
this section, Theorem \ref{2nd.main.thm}, shows that $K[t^{\pm 1};\s]$ contains a free 
subalgebra if and only if $\lim_{n \to \infty}\rho((\s^n)^*)^{1/n} >1$ where 
$\s^*:\NS(X)_\RR \to \NS(X)_\RR$ is the map induced by a birational map $\s:X \dashrightarrow X$ that induces 
$\s \in \Aut_k(K)$.

The proof of Theorem \ref{2nd.main.thm} is similar to that of Theorem \ref{thm.main} but relies in an essential way on 
results of Diller and Favre \cite{DF}, and Rogalski's extension of their results to arbitrary algebraically closed fields \cite{Rog}. 
It is also necessary to replace the the finite dimensional vector space $\NS(X)_\RR$ by an infinite dimensional Hilbert space 
that is constructed from the direct limit of the N\'eron-Severi groups for all models for $K$. 

\subsection{Dynamical degree}
\label{sect.dyn.deg}
 
  The {\sf dynamical degree} of  a  birational map $\s:X \dashrightarrow X$ is the number in (\ref{eq.dyn.deg}). We denote it by 
  $\l(\s)$.

\begin{prop}
\cite[Prop. 3.1]{BFJ}
\label{prop.BFJ}
Let $\s:X \dashrightarrow X$ be a birational map. If  $||\cdot ||$ is a matrix norm on  
$\End_\CC\big(\NS(X)_\CC\big)$, then  
\begin{equation}
\label{eq.dyn.deg}
 \lim_{n \to \infty}\sup  ||(\s^n)^*||^{1/n} \; = \; \lim_{n \to \infty}\sup \big( (\s^n)^*[H] \cdot [H] \big)^{1/n}.
\end{equation}
for every ample divisor $H$ on $X$.
\end{prop}

If $k=\CC$ and $\s:X \to X$ is a birational morphism, then $\l(\s)=\rho(\s^*)$ so $\l(\s)>1$ if and only if $\s$ has positive entropy.

 \subsection{Non-geometric automorphisms}
A {\sf model} for $K$ is a smooth projective surface $X$ such that $K\cong k(X)$.   
Following Rogalski \cite[p.5922]{Rog}, we call a $k$-algebra automorphism $\s \in \Aut_k(K)$ {\sf geometric}
if it is induced by an automorphism  $\s'$ of a model for $K$; we then call $\s'$ a {\sf geometric model} for $\s$. 
If $\s'$ is a geometric model for $\s$, then Theorem \ref{2nd.main.thm} follows from Theorem \ref{main.thm}
because $\lim_{n \to \infty}\rho((\s^n)^*)^{1/n} = \rho(\s'^*)$.

The following result is known but we include a proof for the convenience of the reader.

\begin{prop}
Let $k$ be an algebraically closed field, $K/k$ a finitely generated extension of transcendence degree 2, and $\s$ a non-geometric automorphism of $K/k$. Either $K \cong k(\PP^2)$ or $K \cong k(C \times \PP^1)$ where $C$ is a 
smooth curve of genus $\ge 1$.  If $\l(\s)>1$, then $K \cong k(\PP^2)$. 
\end{prop}
\begin{pf}
Let $X$ be a minimal model for $K$ and $\s:X \dashrightarrow X$ a birational map that
induces the given non-geometric automorphism of $K/k$.
Let $K_X$ be the canonical class for $X$. If $K_X$ were nef, then 
$\s$ would be an automorphism by \cite[Thm. 2, p.180]{IS}. Hence $K_X$ is not nef.
Therefore, by  \cite[Thm.1, p.176]{IS}, $X$ is either $\PP^2$ or a ruled surface, i.e., a $\PP^1$-bundle over a curve. 
Hence $K$ is isomorphic to $k(\PP^2)$ or $k(C \times \PP^1)$ where $C$ is a smooth curve of genus $\ge 1$. 

Let $C$ be a smooth curve of genus $\ge 1$. As stated in the last paragraph of the proof of \cite[Lem. 4.2]{DF}, 
every birational map $C \times \PP^1 \dashrightarrow C \times \PP^1$ has dynamical degree 1. Thus, $k(\PP^2)$ is the only 
field of transcendence degree two having a non-geometric automorphism of dynamical degree $>1$.  
\end{pf}

\subsubsection{The degree of a birational map $\PP^2 \dashrightarrow \PP^2$}
If $\s:\PP^2 \dashrightarrow \PP^2$ is a birational map, there are homogeneous polynomials $f_1,f_2,f_3 \in k[x,y,z]$ having the 
same degree such that $\gcd\{f_1,f_2,f_3\}=1$ and $\s(p)=\big(f_1(p),f_2(p),f_3(p)\big)$ for all but finitely many $p \in \PP^2$. 
The 
{\sf degree} of $\s$ is the common degree of $f_1,f_2,f_3$.

\subsubsection{}
 \label{ssect.non-geom}
Almost every birational map $\PP^2 \dashrightarrow \PP^2$ of degree $\ge 2$ has
dynamical degree $>1$. Precisely, the set of birational maps $\PP^2 \dashrightarrow \PP^2$ of degree $d \ge 2$
having dynamical degree $>1$ is Zariski dense in the set of  all birational maps $\PP^2 \dashrightarrow \PP^2$ of degree $d$
 \cite[Thm.1.6]{X}.

\subsubsection{}
\label{ssect.Henon.1} 
In \S\ref{ssect.Henon.2}, we show that the  automorphism of $\CC(x,y)$
induced by the H\'enon map $\CC^2 \to \CC^2$, $(x,y) \mapsto (1+y-ax^2,bx)$, where $a\ne 0$ and $b$ are 
fixed complex numbers, is not geometric.

 Artin and Rogalski,   \cite[Ex. 3.6]{Rog}, exhibit a birational map $\PP^2 \dashrightarrow \PP^2$ that is not conjugate to an automorphism, namely
$\phi\s:\PP^2 \dashrightarrow \PP^2$ where $\phi$ is a generic element of $\PGL(2,k)$ and $\s$ is the 
quadratic Cremona map $(x,y,z) \mapsto \big(\frac{1}{x},\frac{1}{y},\frac{1}{z}\big)$. 

More non-geometric birational maps  $\PP^2 \dashrightarrow \PP^2$ appear at \cite[Prop. 9.3]{DF}. 

Every ruled surface has a non-geometric  birational self-map that preserves the ruling
\cite[Rmk. 7.3]{DF}.

 \subsection{}
Consider a smooth projective surface $X$ endowed with a birational map $\s: X \dashrightarrow X$. 
By blowing up the finite set of points where $\s$ is not defined one can obtain 
birational morphisms $f$ and $g$ that fit into a commutative diagram
$$
  \UseComputerModernTips
\xymatrix{
& Y \ar[dl]_{f} \ar[dr]^{g}
\\
X \ar@{-->}[rr]_{\s} && X.
}
$$
The group homomorphisms 
$\s_*,\s^*:\NS(X)_\RR \to \NS(X)_\RR$ defined by $\s_*:=g_*f^*$ and $\s^*:=f_*g^*$ do not depend on the choice of $f$, $g$, or $Y$. The maps $\s^*$ and $\s_*$ will not always preserve the intersection form but they are adjoint to each other
in the sense that $(\s^*u)\cdot v = u\cdot (\s_*v)$.

When $\s$ is not a morphism $(\s^*)^n$ may not equal $(\s^n)^*$ so some of the arguments
in \S\ref{sect.proof} can not be used.

\subsubsection{}
If $\sigma:\PP^2 \dashrightarrow \PP^2$ is the quadratic Cremona map $\sigma(x,y,z):= \big(\frac{1}{x},\frac{1}{y},\frac{1}{z}\big)$, then
$(\sigma^*)^2 \ne (\sigma^2)^*$. The action of $\sigma^*$ on $\Pic(\PP^2)=\NS(\PP^2)=\ZZ$  is multiplication by $2$ 
but $(\sigma^2)^*$ is multiplication by 1 because $\sigma^2=\id_{\PP^2}$.
Although $\rho(\sigma^*)>1$,  $k(\PP^2)[t^{\pm 1};\sigma]$ does not contain a free subalgebra
because $k(\PP^2)[t^{\pm 1};\sigma]$ is a free module of rank 2 over its commutative subalgebra $k(\PP^2)[t^{\pm 2}]$.


 Theorem \ref{thm.DF.0.1} and Lemma \ref{lem.DF.Rog} imply that 
$$
\l(\s) \, = \,  \lim_{n \to \infty}\rho((\s^n)^*)^{1/n} \,  = \, \lim_{n \to \infty}\rho((\s^n)_*)^{1/n}.
$$

\subsubsection{}

The dynamical degree of the quadratic Cremona map is 1.  
This suggests that it is $\l(\s)$, not $\rho(\s^*)$, that
determines whether $k(X)[t^{\pm 1}; \s]$ contains a free subalgebra.

 \subsection{Stable maps and conjugate maps}
The next three results are due to Diller and Favre \cite{DF} when $k=\CC$, and to Rogalski for arbitrary algebraically closed fields \cite{Rog}. 

A birational map $\tau: Y \dashrightarrow Y$ is {\sf stable} if it has the properties in the next theorem. By \cite[Rmk. 2.7]{Rog}, 
$\tau^{-1}$ is stable if $\tau$ is. 
 
\begin{thm}
\cite[Thm. 1.14]{DF} \cite[Lem. 2.8]{Rog}
\label{thm.DFR}
Let $Y$ be a smooth projective surface and $\tau: Y \dashrightarrow Y$ a birational map. 
The following are equivalent:
\begin{enumerate}
  \item 
  $\tau^n(C)$ is not a fundamental point of $\tau$ for any irreducible curve $C \subset Y$ or any integer $n \ge 1$;
  \item 
  $(\tau^*)^n=(\tau^n)^*$  for all $n\in \ZZ$.
\end{enumerate}
If these conditions hold, then $(\tau_*)^n=(\tau^n)_*$ for all $n \ge 1$.
\end{thm}
 
If $\s: X \dashrightarrow X$ is a birational map and  $\pi: Y \dashrightarrow X$ is another birational map,  
we say that the birational map $  \pi^{-1}\s\pi: Y \dashrightarrow Y$ is {\sf conjugate} to $\s$.

\begin{thm}
\label{thm.DF.0.1}
\cite[Thm. 0.1]{DF} \cite[Thm. 2.9]{Rog} 
Every birational map $\s: X \dashrightarrow X$ is conjugate to a stable map.
\end{thm}

\begin{lem}
\cite[Cor. 1.16]{DF}
\cite[Lem. 2.11]{Rog}
\label{lem.DF.Rog}
\begin{enumerate}
  \item 
  If $\tau:Y \dashrightarrow Y$ is conjugate to $\s:X \dashrightarrow X$, then $\l(\tau)=\l(\s)$.
  \item 
  If  $\tau:Y \dashrightarrow Y$ is stable, then $\l(\tau)=\rho(\tau^*)$.
\end{enumerate}
\end{lem}

\subsubsection{}
\label{ssect.Henon.2} 
I am grateful to Eric Bedford for the following proof of the fact that the automorphism of $\CC(x,y)$
induced by the H\'enon map  is not geometric. 
By \cite[Thm. 2.1]{FM}, the dynamical degree of the birational map  $\s: \PP^2 \dashrightarrow \PP^2$ induced by 
$(x,y) \mapsto (1+y-ax^2,bx)$ is 2. If $\s$ were conjugate to an automorphism $\tau:Y \to Y$, then 
$\rho(\tau^*)=\l(\tau)=\l(\s)=2$. However, $\tau^*$ acts on $\NS(Y)_\RR$ preserving a lattice so its
characteristic polynomial is a monic polynomial with integer coefficients. The eigenvalues of $\tau^*$ are therefore
algebraic integers. By Proposition \ref{prop.q.form},  $\rho(\tau^*)$ is an eigenvalue of $\tau^*$ 
so is certainly not equal to 2.

\subsection{}
We will use the fact that the isomorphism 
class of the graded algebra $k(X)[t^{\pm 1};\s]$ is not changed when $\s$ is replaced by a map conjugate to it. 

The following observation is well-known and easy to verify. 

\begin{lem}
\label{lem.isom}
Let $\pi$ and $\s$ be automorphisms of a field $K$. If $\tau=\pi^{-1}\s\pi$, then the map $\Psi:K[t^{\pm 1};\s] \to 
K[t^{\pm 1};\tau]$ defined by $\Psi(at^j)=\pi^{-1}(a)t^j$ for all $a \in K$  is an isomorphism of 
graded $k$-algebras.
\end{lem}

 \subsection{The Picard-Manin group}

When $\s:X \dashrightarrow X$ is not geometric the argument in Theorem \ref{main.thm} 
must be adapted by replacing $\NS(X)$ and $\NS(X)_\RR$ by suitable direct limits of $\NS(X')$ 
and $\NS(X')_\RR$ taken over the directed system of all homomorphisms $\ve^*:\NS(X) \to \NS(X')$ induced by
all birational morphisms $\ve:X' \to X$.

The {\sf Picard-Manin} group of $X$ is 
$$
{\sf Z}(X):= \liminj \NS(X').
$$
If $\ve:X' \to X$ is a birational morphism the map $\ve^*:\NS(X) \to \NS(X')$ is injective so ${\sf Z}(X)$ is 
the directed union of its subgroups $\NS(X')$. Furthermore,   
$(\ve^* d) \cdot (\ve^* d')= d \cdot d'$ for all $d,d' \in \NS(X)$ so there is an induced intersection form on  ${\sf Z}(X)$
that restricts to the usual intersection form on each $\NS(X')$. 

The intersection form extends to ${\sf Z}(X)_\RR:= {\sf Z}(X) \otimes_\ZZ \RR$ and there is an associated 
non-degenerate quadratic form on ${\sf Z}(X)_\RR$. If $\omega \in {\sf Z}(X)_\RR$ is such that $\omega \cdot \omega >0$,
then the quadratic form is negative definite on its orthogonal complement, $\omega^\perp$.

The completion of ${\sf Z}(X)_\RR$ with respect to the metric induced by the quadratic form is denoted by $\cZ(X)$. 
As remarked at \cite[p.525]{BFJ}, if $\omega \in {\sf Z}(X)$ has the property that $\omega\cdot \omega >0$, then 
$\cZ(X)=\RR \omega \oplus \overline{\omega^\perp}$ and $\cZ$ becomes a 
Hilbert space with respect to the norm $||t\omega \oplus \a ||^2=t^2-(\a\cdot \a)$. The norm depends on the choice of $\omega$
but the structure of $\cZ(X)$ as a topological vector space does not depend on the choice of $\omega$.

\subsubsection{}
The construction of ${\sf Z}(X)$ first appeared in \S34 of Manin's book
{\it Cubic Forms} \cite{Man}. Other accounts, in the context of complex dynamics, appear in Cantat's paper \cite{Can11}
and, using the notation $C({\frak X})$ for ${\sf Z}(X)_\RR$, and $L^2({\frak X})$ for $\cZ(X)$, in \cite{BFJ}. 
We refer the reader to to those sources for more detail and proofs.

\subsubsection{}
If $\pi:Y \to X$ is a birational morphism, then the directed systems used to define ${\sf Z}(Y)$ and ${\sf Z}(X)$ are cofinal
so there is an induced isomorphism $\pi_*:{\sf Z}(Y) \to {\sf Z}(X)$. This isomorphism has the property that $(\pi_*d) \cdot
(\pi_*d')=d\cdot d'$.  Thus ${\sf Z}(X)$ is an invariant of $k(X)$. 

\subsubsection{}
The {\sf Riemann-Zariski} surface \cite[Ch. V, \S17]{ZS} is the inverse limit
$$
{\frak X}:=\liminj X_\ve
$$
taken over the directed system of all birational morphisms $\ve:X_\ve \to X$. Intuitively, ${\frak X}$ is obtained by blowing up all
points on all surfaces birationally isomorphic to $X$. In a sense, ${\cZ}(X)$ plays the role of the N\'eron-Severi group of ${\frak X}$.

\subsubsection{}
Let ${\sf Bir}(X)$ denote the group of birational maps $X \dashrightarrow X$. Thus, ${\sf Bir}(X)$ is the group of $k$-algebra
automorphisms of $k(X)$. 

Let $\s \in {\sf Bir}(X)$. By blowing up the indeterminacy locus of $\s$ one obtains a commutative diagram
$$
  \UseComputerModernTips
\xymatrix{
& Y \ar[dl]_{f} \ar[dr]^{g}
\\
X \ar@{-->}[rr]_{\s} && X.
}
$$
in which $f$ and $g$ are birational morphisms. The morphisms $f$ and $g$ induce isomorphisms
$f_*,g_*:{\sf Z}(Y) \to {\sf Z}(X)$ so we may define 
$$
\s_*:=(g_*) \circ (f_*)^{-1} \, \in \,  \GL({\sf Z}(X)).
$$

Consider a smooth projective surface $X$ endowed with a birational map $\s: X \dashrightarrow X$. 
By blowing up the finite set of points where $\s$ is not defined one can obtain 
birational morphisms $f$ and $g$ that fit into a commutative diagram

The group homomorphisms 
$\s_*,\s^*:\NS(X)_\RR \to \NS(X)_\RR$ defined by $\s_*:=g_*f^*$ and $\s^*:=f_*g^*$ do not depend on the choice of $f$, $g$, or $Y$. The maps $\s^*$ and $\s_*$ will not always preserve the intersection form but they are adjoint to each other
in the sense that $(\s^*u)\cdot v = u\cdot (\s_*v)$.

\begin{thm}
\cite[Thm. 34.8]{Man}
The map $\s \mapsto \s_*$ is an injective group homomorphism ${\sf Bir}(X) \to \GL({\sf Z}(X))$. Furthermore, $\s_*$
preserves the intersection form on ${\sf Z}(X)$ and the nef cone.
\end{thm}

In particular,   $(\s^n)_* = (\s_*)^n:{\sf Z}(X) \to {\sf Z}(X)$ for all $n \in \ZZ$.

The functoriality argument in \cite[\S2]{BFJ} shows that the map $\s \mapsto \s_*$ extends to a group homomorphism from 
${\sf Bir}(X)$ to the isometry group of $\cZ(X)$.

\begin{thm}
\cite[Thm. 3.2]{BFJ}
\label{thm.BFJ}
If $\s: X \dashrightarrow X$ is a birational map with dynamical degree $\l$, then there are non-zero nef classes 
$e_{{}_+}, e_{{}_-} \in \cZ(X)$   such that $\s^*(e_{{}_+})=\l e_{{}_+}$ and 
$\s^*(e_{{}_-})=\l^{-1} e_{{}_-}$. 
\end{thm}

\subsubsection{}
We make no use of the nef cone in $\cZ(X)$ and only mention it for the sake of
completeness.

\subsection{}
The following facts, which will be used in the proof of Theorem \ref{2nd.main.thm},  
are simple consequences of Bertini's Theorem \cite[8.20.2]{Ha} and the uncountability of $k$.

\begin{prop}
\label{prop.avoiding.S}
Let $Y$ be a smooth projective surface over an uncountable algebraically closed field. 
Let $\tau:Y \dashrightarrow Y$ be a 
birational map. Let
$$
S=\{y \in Y \; | \; \hbox{some $\tau^n$, $n \in \ZZ$, is not defined at $y$}\}.
$$
\begin{enumerate}
  \item 
the set $S$ is countable;  
\item
the number of curves that are contracted to a point by some $\tau^n$, $n \in \ZZ$, is countable;
\item 
there are uncountably many very ample effective divisors $H$ such that $H \cap S=\varnothing$. 
\item 
Let $H$ be a very ample effective divisor such that $H \cap S=\varnothing$ and let $h \in k(Y)$ be a rational function whose divisor of zeroes is $H$. Then  there are uncountably many smooth irreducible curves $C \in |H|$ such that 
\begin{enumerate}
  \item 
  $C \cap S=\varnothing$ and
  \item
  $C$ is not contracted to a point by any $\tau^n$, $n \in \ZZ$, and
  \item 
 $C$ is not contained in the zero locus  of $\tau^n(h)$ or the zero locus  of $\tau^n(h^{-1})$ for any $n \in \ZZ$.
\end{enumerate}
\end{enumerate}
\end{prop}

\begin{thm}
\label{2nd.main.thm}
Let $X$ be a smooth projective surface over an uncountable algebraically closed field $k$. 
Let $\s: X \dashrightarrow X$ be a birational map and $\tau:Y \dashrightarrow Y$ a stable birational map conjugate to $\s$.
If $\l(\s)>1$,  then $k(X)[t^{\pm 1};\s]$ contains a free subalgebra generated by two homogeneous elements of degree $n$ 
for every integer $n$ such that $\l(\s)^n \ge 5+2\sqrt{6}$. More explicitly,
\begin{enumerate}
\item
$k(X)[t^{\pm 1};\s] \cong k(Y)[t^{\pm 1};\tau]$;
\item
if $H$ is a very ample effective divisor on having the properties in Proposition \ref{prop.avoiding.S} and $h \in k(Y)$ is such that $(h)_0=H$, then 
$k\{t^n,ht^n\}$ is a free subalgebra of $k(Y)[t^{\pm 1};\tau]$;
  \item 
  $k\{at^n,bt^n\}$ is a free subalgebra of $k(Y)[t^{\pm 1};\tau]$ if $a,b \in k(Y)$ are such that $a^{-1} b=h$. 
\end{enumerate}
\end{thm}
\begin{pf}
First, replace $\s$ by a power of itself such that its dynamical degree is $\l \ge 5+2\sqrt{6}$. 
Let $\tau:Y \dashrightarrow Y$ be a stable birational map conjugate to $\s$.  
We have $\l=\l(\s)=\l(\tau)=\rho(\tau^*)$ where $\tau^*:\cZ(Y) \to \cZ(Y)$. By Lemma \ref{lem.isom},  
$k(X)[t^{\pm 1};\s] \cong k(Y)[t^{\pm 1};\tau]$.

Let $H$ be a very ample effective divisor and $C$ a smooth irreducible curve on $Y$ that have the properties listed in
Proposition \ref{prop.avoiding.S}.
 
We will now explain how the proof of Theorem \ref{main.thm} can be modified to show that
\begin{equation}
\label{H.C}
 ( \tau^*)^{j+1}[H] \! \cdot \!  [C] \, > \,  2 (\tau^*)^{j}[H] \! \cdot \!  [C] 
 \end{equation}
  for all $j \ge 0$. 

By Theorem \ref{thm.BFJ}, there are non-zero elements $e_{{}_+},e_{{}_-} \in \cZ(Y)$ 
such that $\tau^*(e_\pm)=\l^{\pm 1}e_{\pm}$. Because $\l>1$ and $(\tau^* e_{{}_{\pm}}) \cdot (\tau^* e_{{}_{\pm}}) = 
e_{{}_{\pm}} \cdot  e_{{}_{\pm}}$, $ e_{{}_+}\!\! \cdot e_{{}_+} = e_{{}_-}\! \! \cdot e_{{}_-}=0$. 
The signature of the quadratic form on $\cZ(X)$ is $(1,\infty)$ so $\RR e_{{}_+} \oplus \RR e_{{}_-}$ is not totally isotropic. 
Therefore $e_{{}_+}\!\! \cdot  e_{{}_-} \ne 0$.  Let $\cW=(\RR e_{{}_+} \oplus \RR e_{{}_-})^\perp$ denote the orthogonal complement of $\RR e_{{}_+} \oplus \RR e_{{}_-}$ in $\cZ(Y)$. 
Because the quadratic form on $\cZ(Y)$  has signature $(1,\infty)$,    $w\! \cdot \!  w\le 0$ for all $w \in \cW$. 
Since $\tau^*$ is orthogonal,  $\tau^*(\cW)=\cW$.

Write $[H]= \a e_{{}_+} + \b e_{{}_-} + w$ for its class in $\cZ(Y)$ where 
$w \in \cW$. Because $H\!\! \cdot \!  H >0$, $\a\b(e_{{}_+}\!\! \cdot \!  e_{{}_-})> 0$. Replace $e_{{}_+}$ and $e_{{}_-}$ by multiples of themselves so that $[H]= e_{{}_+} \, + \, e_{{}_-} \, +\, w$. Now one follows the argument in Theorem \ref{main.thm} to 
prove (\ref{H.C}); the only change is that one must use the Cauchy-Schwarz inequality in $\cZ(Y)$.  

The argument in the proof of Theorem \ref{main.thm} involving the passage to the 
  discrete valuation ring $\cO_{Y,C}$ and then to its residue field $k(C)$ can be used to complete 
  the proof of this theorem.  
\end{pf}

\begin{prop}
Let $X$ be a smooth projective surface over an algebraically closed field $k$ and 
$\s:X \dashrightarrow X$ a birational map. If $\l(\s)=1$, then $k(X)[t^{\pm 1};\s]$ does not contain a free 
subalgebra on $\ge 2$ variables.
\end{prop}
\begin{pf}
By \cite[Thm. 4.1]{DF} in the case $k=\CC$ and by \cite[Prop. 2.14]{Rog} for an arbitrary algebraically closed field $k$, $\s$ 
is conjugate to an automorphism $\tau:Y \to Y$ such that $(\tau^n)^*$ is the identity. 
Since $1=\l(\s)=\l(\tau)=\rho(\tau^*)$,  Theorem \ref{thm.RZ} tells us that $k(X)[t^{\pm 1};\tau]$ 
does not contain  a free subalgebra.  This completes the proof because $k(X)[t^{\pm 1};\tau] \cong k(X)[t^{\pm 1};\s]$.
\end{pf}

\section{Examples}
\label{sect.egs}

\subsection{}
Compact complex surfaces having automorphisms of positive entropy are rather rare.

Suppose, for example, that $v \in \NS(X)_\RR$ is an eigenvector for $\s^*$ having eigenvalue $\l>1$. 
The canonical class, $[K_X]$, in $\NS(X)_\RR$ is fixed by $\s^*$, so $\l [K_X] \cdot v = \s^*[K_X]\cdot v = [K_X] \cdot v$
which implies that $[K_X] \cdot v=0$. Hence $K_X$ is neither ample nor anti-ample. 

\begin{prop}
[Cantat]
\cite[Prop.1]{Can99} \cite[\S2.3]{Can2}
\label{prop.cantat}
Let $X$ be a compact complex surface having an automorphism of positive topological entropy. After contracting any exceptional curves that have finite order under the automorphism, $X$ is either an abelian surface, a K3 surface, an Enriques surface, 
or  a rational surface obtained by blowing up at least 10 points of $\PP^2$.
\end{prop}

The Kodaira dimension of a rational surface is $-\infty$, and the other surfaces in Proposition
\ref{prop.cantat} have Kodaira dimension 0.

\subsection{}


It is well-known that $\CC\PP^2$ has no automorphisms of positive entropy. 
However, after blowing up $\CC\PP^2$ at a suitable set of $\ge 10$ points one arrives at a rational surface 
that does have an automorphism of  positive entropy---see, for example,  Bedford and Kim  \cite{BK} and McMullen \cite{M2}.

\subsection{An abelian surface}
\label{sect.EE}

Let $\omega=e^{2\pi i/3}$, $\L=\ZZ[\omega]$, $E=\CC/\L$, and $X=E \times E$. By McMullen  \cite[Thm. 1.3]{M3},
the automorphism of $\CC^2$ given by multiplication by 
$$
\begin{pmatrix}
  1    &  1  \\
   \omega   & 0 
\end{pmatrix}
$$
induces an automorphism of $X$ whose entropy is $\log(\l)$ where $\l$ is the largest root of $x^4-x^3-x^2-x+1$. 
Here $\l \approx 1.72208380$. The curve $E$ can be defined by the equation $y^2z=x^3+z^3$ so 
$k(X)=k(x_1,x_2,y_1,y_2)$ with relations $y_i^2=x_i^3+1$. The divisor of zeroes of $x_1x_2$ is very ample so
$k\{t^n,x_1x_2t^n\}$ is a free subalgebra of $k(X)[t^{\pm 1};\s]$.

\subsection{K3 surfaces with Picard number 2}
\label{sect.pic=2}

We now consider a  family of K3 surfaces $S$ first examined by Wehler \cite{W} and further studied by Baragar, \cite{B2} and \cite{B3}, Silverman \cite{S}, and others.
Each $S$ is the common zero locus of a generic $(1,1)$- and a generic $(2,2)$-form on
 $\PP^2\times \PP^2$.  Proofs of the facts we need about $S$ can be found in \cite{W} and \cite{S}.

 The divisors  $\{H_1,H_2\}$ that are the zero loci of a $(1,0)$- and a $(0,1)$-form on $\PP^2\times \PP^2$ are a basis for
  $\Pic(S) =\NS(S) \cong \ZZ^2$.
  
Let $p_i:S \to \PP^2$, $i=1,2$, be the projections $p_i(x_1,x_2)=x_i$. These projections have degree two so induce
involutions of $S$ that interchange the two ``sheets'' of the covering. In more detail, because of the $(2,2)$-form 
vanishing on $S$ every fiber of $p_i^{-1}$ consists of two points counted with multiplicity. Let $\s_1:S \to S$ be the involution 
defined by $\s_1(x_1,x_2)=(x_1,x_2')$ where $p_1^{-1}(x_1)=\{(x_1,x_2),(x_1,x_2')\}$. There is a similar involution $\s_2$.
 By \cite[Thm. 2.9]{W},  when $S$ is generic, $\s_1$ and $\s_2$ generate $\Aut(S)$ which is isomorphic to 
 the free product $\ZZ_2 * \ZZ_2$.
 
Let  $\s=\s_2\s_1$.
 
 \begin{prop}
 \label{prop.eg.pic=2}
 Let $h \in k(S)$ be such that $(h)_0=H_1+H_2$. If $f,g \in k(S)$ are such that $f=gh$, then $k\{ft,gt\}$ is a free subalgebra
 of $k(S)[t^{\pm 1};\s]$.
 \end{prop}
 
 The proof of Proposition \ref{prop.eg.pic=2} occupies the rest of this section. 
It rests on Theorem \ref{main.thm} and some calculations in $\NS(S)_\RR$
that appear in \cite{S} and \cite{W}.

 \subsubsection{}
 With respect to the basis $\{H_1,H_2\}$, the action of $\s^*$ on $\NS(S)_\RR$ is given by the matrix 
${{15 \; \; 4 } \choose {-4 \; \; 1}}$. The eigenvalues of $\s^*$ are $7 \pm 4 \sqrt{3}$. 
 
 \begin{lem}
 \cite[Lemma 2.1]{S}
 Let $\a = 2 + \sqrt{3}$ and define $e_{{}_+}=\a H_1-H_2$ and $e_{{}_-}= -H_1+\a H_2$. Then
 $$
 \begin{array}{ll}
 \s_1^*H_1 =  H_1  \phantom{xxxx} &  \s_1^*H_2 =  4H_1 - H_2  \\
  \s_2^*H_1 =  -H_1 +4H_2  \phantom{xxxx} &  \s_2^*H_2 =  H_2  \\ 
  \s_1^*e_{{}_+} =  \a^{-1}e_{{}_-}   \phantom{xxxx}  &  \s_1^*e_{{}_-} =  \a e_{{}_+}  \\
   \s_2^*e_{{}_+} = \a  e_{{}_-}   \phantom{xxxxxx}  &   \s_2^*e_{{}_-} =  \a^{-1}e_{{}_+} \\
   \s^*e_{{}_+}=\a^2 e_{{}_+}  \phantom{xxxxxx}  &  \s^*e_{{}_-}=\a^{-2} e_{{}_-}.
 \end{array}
 $$
 \end{lem}  
 
 The intersection matrix  with respect to  $\{H_1,H_2\}$  is ${{2 \; \; 4 } \choose {4 \; \; 2}}$ so
 $$
 \begin{array}{ll}
 H_1 \cdot H_1 = H_2 \cdot H_2 = 2  \phantom{xxxx} &  H_1 \cdot H_2=4  \\
  e_{{}_+} \cdot e_{{}_+}=e_{{}_-}\cdot e_{{}_-} = 0   \phantom{xxxx}  &  
  e_{{}_+}\cdot  e_{{}_-} = 12\a \\ 
  e_{{}_+}\cdot H_1 =    e_{{}_-}\cdot H_2  &  e_{{}_+}\cdot H_2 =    e_{{}_-}\cdot H_1
  \\
\phantom{  e_{{}_+}\cdot H_1} =2 \sqrt{3} = \a+\a^{-1} \phantom{xxxx} &\phantom{  e_{{}_+}\cdot H_1}   = 6+4\sqrt{3}    =\a^2-1.
 \end{array}
 $$

\begin{prop}
\label{prop.good.curve}
The divisor $H:=H_1+H_2$ is very ample and there is a smooth irreducible curve $C$ 
such that $\s^{n+1}H \cdot C > 2\,\s^nH \cdot C$ for all $n \ge 0$.
\end{prop}
\begin{pf}
Since $H_1+H_2$ is the intersection of $S$ with a very ample divisor on $\PP^2 \times \PP^2$, $H$ is a very ample divisor on 
$S$. Bertini's hyperplane theorem implies that the linear equivalence class of $H$ contains a 
smooth irreducible curve; let $C$ be such a curve. Let $\l=\a^2$. Since $e_{{}_+} + e_{{}_-} =\a H$, 
\begin{align*}
(\s^*)^nH \cdot C  & = (\s^*)^nH \cdot H
\\
&= \a^{-2}(\s^*)^n(e_{{}_+} + e_{{}_-} )  \cdot (e_{{}_+} + e_{{}_-} )
\\
&= \a^{-2} (\l^{n}e_{{}_+} + \l^{-n} e_{{}_-} )  \cdot (e_{{}_+} + e_{{}_-} )
\\
&= 12 \a^{-1}(\l^{n}+ \l^{-n})
\end{align*}
The inequality in the statement of the proposition is easily verified.  
\end{pf}

\subsubsection{The proof of Proposition \ref{prop.eg.pic=2}}
Let $((x,y,z),(x',y',z'))$ be bihomogeneous coordinate functions on $\PP^2 \times \PP^2$. Since
$H$ is the zero locus of a (1,1)-form,  $p:=(ax+by+cz)(a'x'+b'y'+c'z')$ say, it is the divisor of zeroes of a rational
function $h=p/q$ where $q$ is any (1,1)-form that is relatively prime to $p$. 

It follows from Theorem \ref{main.thm} that $k\{t,ht\}$ is a free subalgebra of $k(S)[t^{\pm 1};\s]$ and Proposition \ref{prop.eg.pic=2}
now follows from Lemma \ref{lem.known}.

\subsection{K3 surfaces with Picard number 3}

Let $X$ be a smooth hypersurface in $\PP^1 \times \PP^1 \times \PP^1$ cut out by a $(2,2,2)$ form. 
The adjunction formula implies that the canonical bundle of $X$ is trivial. 
From the exact sequence $0 \to \cO(-2,-2,-2)\cong \cI  \to \cO \to \cO_X \to 0$ and the fact that $h^1(\cO)=h^2(\cI)=0$ one sees that $H^1(X,\cO_X)=0$.  Hence $X$ is a K3 surface. 

At the end of \cite{W}, Wehler remarks that the automorphism group of the generic member of this family is isomorphic 
to $\ZZ_2 * \ZZ_2 * \ZZ_2$ (generated by the involutions $\s_1$, $\s_2$, $\s_3$, defined in the next paragraph). Such $X$ have been studied in detail by Baragar \cite{B1} and \cite{B3}, Cantat-Oguiso \cite[Sect. 3.4]{CO}, McMullen \cite{M1}, and others. 

A specific example in this family is given by the closure in $\PP^1 \times \PP^1 \times \PP^1$ of the affine variety 
$(1+x^2)(1+y^2)(1+z^2)+Axyz=2$ where $A \in \RR$. 
Some beautiful pictures illustrating the behavior of positive entropy automorphisms on the real surface $X(\RR)$ 
for various values of $A$ appear in Cantat \cite{Can2} and McMullen \cite{M1}. 

\subsubsection{}
Let $p_i:X \to \PP^1 \times \PP^1$, $i=1,2,3$, be the projections obtained by forgetting the $i^{\th}$ coordinate.
Let $H_i=p_i^{-1}(\PP^1 \times \PP^1)$. 
Because $X$ has tri-degree $(2,2,2)$ the preimage under $p_i$ of a point in $\PP^1 \times \PP^1$ consists of two points, 
counted with multiplicity.  It follows that there are three involutions $\s_1$, $\s_2$, $\s_3$, of $X$. For example, 
if $x=(x_1,x_2,x_3) \in X$ and $p_1^{-1}(p_1(x))=\{x,x'=(x_1',x_2,x_3)\}$, then $\s_1(x)=x'$. 

The classes of $H_1$, $H_2$, and $H_3$ form a basis for $\Pic(X)$, $\NS(X)$, and $\NS(X)_\RR$. The intersection form with respect to $\{H_1,H_2,H_3\}$ is given by the matrix
$$
\begin{pmatrix} 0 & 2 & 2 \\ 2 & 0 & 2 \\ 2 & 2 & 0 \end{pmatrix}.
$$
The induced actions of  the $\s_i^*$  on $\NS(X)_\RR$ with respect to the ordered basis $\{H_1,H_2,H_3\}$ are:
$$
\s_{1}^*=\begin{pmatrix} -1 & 0 & 0 \\ 2 & 1 & 0 \\2 & 0 & 1 \end{pmatrix}, \qquad
\s_{2}^*=\begin{pmatrix}  1 & 2 & 0\\  0 & -1 & 0 \\ 0 & 2 & 1 \end{pmatrix}, \qquad
\s_{3}^*=\begin{pmatrix} 1 & 0 & 2 \\ 0 & 1 & 2 \\ 0 & 0 & -1 \end{pmatrix}.
$$
Let $\s=\s_3\s_2\s_1$. 

Let $\eta=\frac{1}{2}(1+\sqrt{5})$, the golden ratio.
The eigenvalues for the action of
$$
\s^*=\s_1^*\s^*_2\s_3^* = \begin{pmatrix} -1 & -2 & -6 \\ 2 & 3 & 10 \\ 2 & 6 & 15 \end{pmatrix}
$$
on $\NS(X)_\RR$ are $\l^{\pm 1}=\eta^{\pm 6} =9 \pm 4\sqrt{5}$ and  $-1$,
with corresponding eigenvectors
$$
e_{{}_+}=\begin{pmatrix} -\eta^{-1} \\ 1 \\  \eta \end{pmatrix}, \quad
e_{{}_-}=\begin{pmatrix}  \eta  \\ -1 \\  -\eta^{-1}   \end{pmatrix}, \quad
e=\begin{pmatrix} 1 \\-3 \\ 1 \end{pmatrix}.
$$
 A straightforward computation gives 
 $$
 e_{{}_-}\cdot e_{{}_-}=e_{{}_+}\cdot e_{{}_+}=0, \quad e_{{}_+}\cdot e_{{}_-}=6, \quad e \cdot e = -20, \quad e \cdot e_{{}_+}=e\cdot e_{{}_-}=0.
$$
 
 The divisor $H:=3(H_1+H_2+H_3)$ is very ample. Its class in $\NS(X)_\RR$ is
 $4e_{{}_+}  +4e_{{}_-} - e $. Therefore $(\s^*)^{n}(H)=4\l^ne_{{}_+} + 4\l^{-n}e_{{}_-} +(-1)^{n+1} e$
and 
$$
(\s^*)^n(H) \cdot H = 96(\l^n+\l^{-n})+20(-1)^{n+1}.
$$
It follows that $(\s^*)^{n+1}(H) \cdot H > 2\,(\s^*)^nH\cdot H$ so, as in \S\ref{sect.pic=2}, one can produce a simple explicit 
element $h \in k(X)$ such that $k\{t,ht\}$ is a free subalgebra of $k(X)[t^{\pm 1};\s]$. We leave the routine details to the interested reader.



\subsection{Non-geometric birational maps $\PP^2 \dashrightarrow \PP^2$}
\label{ssect.PP2}

Let $k$ be an uncountable algebraically closed field.

As mentioned in \S\ref{ssect.non-geom}, $k(\PP^2)$  is the only field of transcendence degree two having a non-geometric automorphism of dynamical degree $>1$. However, by \cite[Thm.1.6]{X}, almost all birational maps 
$\PP^2 \dashrightarrow \PP^2$ of degree $\ge 2$ have dynamical degree $>1$.  

 Let $\s:\PP^2  \dashrightarrow \PP^2$ be a non-geometric birational map such that $\l(\s)>1$.  By Theorem \ref{2nd.main.thm}, $k(\PP^2)[t^{\pm 1};\s]$ contains a free subalgebra.

Let $\s:\PP^2 \dashrightarrow \PP^2$ be the H\'enon map, $(x,y) \mapsto (1+y-ax^2,bx)$, $a \ne 0$. 
Since $\l(\s)=2$ the smallest $n$ such that $\l(s)^n > 5+2 \sqrt{6}$ is $n=4$. The zero locus, $H$, 
of $\frac{x}{y}$ in $\PP^2$ satisfies the conditions in Proposition \ref{prop.avoiding.S} so, by Theorem \ref{2nd.main.thm},
$k\{xt^4,yt^4\}$ is a free subalgebra of $k(x,y)[t^{\pm 1};\s]$.  That is not the best possible result: it is shown in \cite{SPS} that 
$k\{xt,yt\}$ is a free subalgebra of $k(x,y)[t^{\pm 1};\s]$. By the next result, the fact that $k\{xt,yt\}$ is free implies 
that $k\{xt^4,yt^4\}$ is also a free algebra.

\begin{prop}
Let $R$ be a commutative $k$-algebra and $\s \in \Aut_k(R)$. If $a,b \in R$ are such that $k\{at^i,bt^i\}$ is a 
free subalgebra of $R[t^{\pm 1};\s]$ so is $k\{at^{im},bt^{im}\}$ for all $m \ge 1$.
\end{prop}
\begin{pf}
By replacing $\s$ by $\s^i$ we can assume that $i=1$. We do this. Thus, we assume that  $k\{at,bt\}$
is a free subalgebra of $R[t^{\pm 1};\s]$ and will show that $k\{at^{m},bt^{m}\}$ is a free algebra for all $m \ge 1$.

Let $V=ka+kb$. To prove the proposition we must show that 
$$
\dim\big(V\s^m(V)\s^{2m}(V) \ldots \s^{mn}(V)\big) = 2^{n+1}.
$$
The dimension of the vector space
$$
W:= \bigotimes_{\stackrel{i=1}{m\not\;\vert \; i}}^{mn-1} \s^i(V)
$$
is $2^{(m-1)n}$.
Multiplication in $R$ gives a surjective map
\begin{equation}
\label{eq.dimW}
W \otimes V\s^m(V)\s^{2m}(V) \ldots \s^{mn}(V) \longrightarrow V \s(V) \ldots  \s^{mn-1}(V) \s^{mn}(V).
\end{equation}
The dimension of the right-hand side of (\ref{eq.dimW}) is $2^{mn+1}$ so
$$
2^{(m-1)n} \dim\big(V\s^m(V)\s^{2m}(V) \ldots \s^{mn}(V)\big) \ge 2^{mn+1}  .
$$
It follows that 
$\dim\big(V\s^m(V)\s^{2m}(V) \ldots \s^{mn}(V)\big) = 2^{n+1}$. This completes the proof that $k\{at^{m},bt^{m}\}$
is a free algebra.
\end{pf}

Several examples of birational maps $\tau:\PP^2 \dashrightarrow \PP^2$ 
for which $k(\PP^2)[t^{\pm 1};\tau]$ contains a free subalgebra appear in \cite{SPS}.
For example, if $\s$ is the automorphism $\s(x,y)=(x,xy)$ of $k(\PP^2)=k(x,y)$, then the subalgebra $k\{xt,yt\}$  of 
$k(x,y)[t^{\pm 1};\s]$ is not free but $k\{xt^n,yt^n\}$ is free for all $n \ge 2$ \cite[Thm. 3.5]{SPS}.


 \end{document}